\documentclass[preprint,12pt]{elsarticle}

\usepackage{amssymb}
\usepackage{graphicx}
\usepackage{amsmath}
\usepackage{amsfonts}
\usepackage{subfigure}
\usepackage{tikz}
\usepackage{IEEEtrantools}
\usepackage[normalem]{ulem}

\usepackage{tikz} 

\usetikzlibrary{positioning}
\usepackage{pgfplots}

\usepackage{array}
\usepackage{mathtools}

\usepackage{IEEEtrantools}
\IEEEeqnarraydefcol{eqnarrayCscript}{\hfil$\scriptstyle}{$\hfil}

\renewcommand{\S}{S}

\newcommand{\ie}{{\it i.e.},\ }

\newtheorem{theorem}{Theorem}

\newtheorem{example}{Example}

\journal{PEIS}

\begin{document}

\begin{frontmatter}



\title{Performance measures for the two-node queue with finite buffers}


\author[label1]{Yanting Chen}
\address[label1]{College of Mathematics and Econometrics, Hunan University, Changsha, Hunan 410082, P.~R.~China\fnref{label}}
\address[label2]{Stochastic Operations Research, University of Twente, P.O. Box 217, 7500 AE Enschede, The Netherlands\fnref{labe2}}
\ead{yantinchen@hnu.edu.cn}

\author[label2]{{Xinwei Bai}}
\ead{x.bai@utwente.nl}

\author[label2]{Richard J. Boucherie}
\ead{r.j.boucherie@utwente.nl}

\author[label2]{Jasper Goseling}
\ead{j.goseling@utwente.nl}

\address{}

\begin{abstract}
We consider a two-node queue modeled as a two-dimensional random walk. In particular, we consider the case that one or both queues have finite buffers. We develop an approximation scheme based on the Markov reward approach to error bounds in order to bound performance measures of such random walks in terms of a perturbed random walk in which the transitions along the boundaries are different from those in the original model and the invariant measure of the perturbed random walk is of product-form. We then apply this approximation scheme to a tandem queue and some variants of this model, for the case that both buffers are finite. We also apply our approximation scheme to a coupled-queue in which only one of the buffers has finite capacity.

\end{abstract}


\end{frontmatter}

The two-node queue is one of the most extensively studied topics in queueing theory. It can be often modeled as a two-dimensional random walk on the quarter-plane. Hence, it is sufficient to find performance measures of the corresponding two-dimensional random walk if we are interested in the performance of the two-node queue. In this work we analyze the steady-state performance of a two node queue for the particular case that one or both of the queues have finite buffer capacity. Our aim is to develop a general methodology that can be applied to any two-node queue that can be modeled as a two-dimensional random walk on (part of) the quarter-plane. 

A special case of the two-node queue with finite buffers at both queues which has been extensively studied so far, is the tandem queue with finite buffers. An extensive survey of results on this topic is provided in~\cite{balsamo2011queueing, perros1994queueing}. Most of these papers focus on the development of approximations or algorithmic procedures to find steady-state system performances such as throughput and the average number of customers in the system. A popular approach used in such approximations is decomposition, see~\cite{asadathorn1999decomposition, gershwin1987efficient}. The main variations of a two-node queue with finite buffers at both queues are: three or more stations in the tandem queue~\cite{shanthikumar1994bounding}, multiple servers at each station~\cite{van1989simple, van2006performance}, optimal design for allocating finite buffers to the stations~\cite{hillier1995optimal}, general service times~\cite{van1987formal, van2004error}, etc. Numerical results of such approximations often suggest that the proposed approximations are indeed bounds on the specific performance measure, however rigorous proofs are not always available. Moreover, these approximation methods cannot be easily extended to a general method, which determines the steady-state performance measure of any two-node queue with finite buffers at both queues.

Van Dijk et al.~\cite{vandijk88tandem} pioneered in developing error bounds for the system throughput using the product-form modification approach. The method has since been further developed by van Dijk et al.~\cite{van1998bounds,vandijk88perturb} and has been applied to, for instance, Erlang loss networks~\cite{boucherie2009monotonicity}, to networks with breakdowns~\cite{van1988simple}, to queueing networks with non-exponential service~\cite{van2004error} and to wireless communication networks with network coding~\cite{goseling2013energy}. An extensive description and overview of various applications of this method can be found in~\cite{vandijk11inbook}.

A major disadvantage of the error bound method mentioned above is that the verification steps that are required to apply the method can be technically quite complicated. Goseling et al.~\cite{goseling2014linear} developed a general verification technique for random walks in the quarter-plane. This verification technique is based on formulating the application of the error bounds method as solving a linear program. In doing so, it avoids completely the induction proof required in~\cite{vandijk88perturb}. Moreover, instead of only bounding performance measures for specific queueing system, the approximation method developed in~\cite{goseling2014linear} accepts any random walk in the quarter-plane as an input.

The main contribution of the current work is to provide an approximation scheme which can be readily applied to approximate performance measures for any two-node queue in which one or both queues have finite buffer capacity. This is based on modifying the general verification technique developed in~\cite{goseling2014linear} for a two-dimensional random walk on a state space that is finite in one or both dimensions. 

We apply this approximation scheme to a tandem queue with finite buffers at both queues. We show that the error bounds for the blocking probability are improved compared with the error bounds for the blocking probability provided in~\cite{vandijk88tandem}. The method in~\cite{vandijk88tandem} is based on specific model modifications. Apart from this, our approximation scheme is more general in the sense that other interesting performance measures could also be obtained easily. This is an advantage over the methods used in~\cite{van1998bounds, vandijk88tandem, vandijk88perturb} where different model modifications are necessary for different performance measures. Moreover, we show that the error bounds can also be obtained for variations of the tandem queue with finite buffers. In particular, we consider the case that one server speeds-up or slows-down when another server is idle or saturated. 

For a two-node queue with finite buffers at both queues, it is also possible to find the invariant measure by solving a system of linear equations. The complexity solving this system is at least $O(L^{2})$, where $L$ is the size of the smallest buffer. We will demonstrate that the approach that is presented in this paper has a complexity that is constant in $L$. This makes it an interesting alternative to solving for the invariant measure by brute force if $L$ is large.
 
Finally, we apply this approximation scheme to a two-node queue with finite buffers at only one queue. In particular, we apply our results to the coupled-queue~\cite{fayolle1979two}. Contrary to~\cite{fayolle1979two}, we consider the case that one of the queues has finite buffer capacity. The numerical results illustrate that our approximation scheme achieves tight bounds.

There are other means to analyze the two-node queue with finite buffers at only or both queues. In particular, the models considered in this paper are instances of quasi-birth-and-death process (QBD) processes and, therefore,  amendable for a solution using the matrix-geometric approach~\cite{Latouche1999Matrix, Neuts1981Matrix}.
There are many variations on the matrix geometric method, in particular in how to compute the rate matrix. However, all methods share a common complexity of $O(L^3)$, where $L$ is the number of phases, which in our case corresponds to the size of the smallest buffer. Therefore, our approach, with constant complexity in $L$, provides a promising alternative to the matrix geometric method for large $L$. A drawback of our approach is that it in general does not give an exact result, but only bounds.

Another important advantage of our work is that it is possible, though outside the scope of the current paper, to extend our approach to queueing networks with more than two queues and more complicated interactions. Such an extension is not possible for the matrix-geometric method. This paper provides the necessary intermediate step in building up our approach from the first ideas in~\cite{goseling2014linear} towards a completely general method that can be applied to queueing networks for which currently no methods exist by which we can analyze them.

The remainder of this paper proceeds as follows. In Section~\ref{sec:model4E}, we present the model and formulate the research problem. In Section~\ref{sec:as4E}, we provide an approximation scheme to bound performance measures for any two-node queue with finite buffers at both queues. We bound performance measures for a tandem queue with finite buffers and some variants of this model in Section~\ref{sec:errorboundsTandemE}. In Section~\ref{sec:onequeue}, we extend the approximation scheme to any two-node queue with finite buffers at only one queue. In Section~\ref{sec:errorboundsSharing4E}, this extended approximation scheme has been applied to a coupled-queue with processor sharing and finite buffers at only one queue. Finally, we provide concluding remarks in Section~\ref{sec:conclusionE}.

\section{Two-node queue with finite buffers at both queues} \label{sec:model4E}

\subsection{Two-node queue with finite buffers at both queues}

The two-node queue with finite buffers at both queues is a queueing system with two servers, each of them having finite storage capacity. If a job arrives at a server which does not have any more storage capacity, then the job is lost. In general, the two queues influence each other, \ie the service rate at one of the queues depends on the number of jobs at the other.
 
Such a queueing system is naturally modeled as a two-dimensional finite random walk, which we introduce next. The connection between the continuous-time queueing system and the discrete-time random walk, obtained through uniformization, is made explicit for various examples in Section~\ref{sec:errorboundsTandemE} and Section~\ref{sec:errorboundsSharing4E}.

\subsection{Two-dimensional finite random walk on both axis}

\begin{figure}
\begin{center}
{
\includegraphics{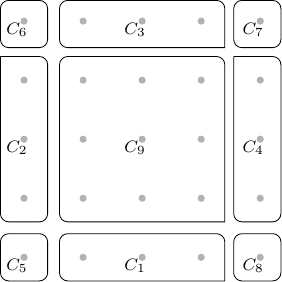}
}
\caption{$C$-partition of $S$ with components $C_1, C_2, \cdots, C_{9}$. \label{fig:partXE}}
\end{center}
\end{figure}

\begin{figure}
\hfill
\includegraphics{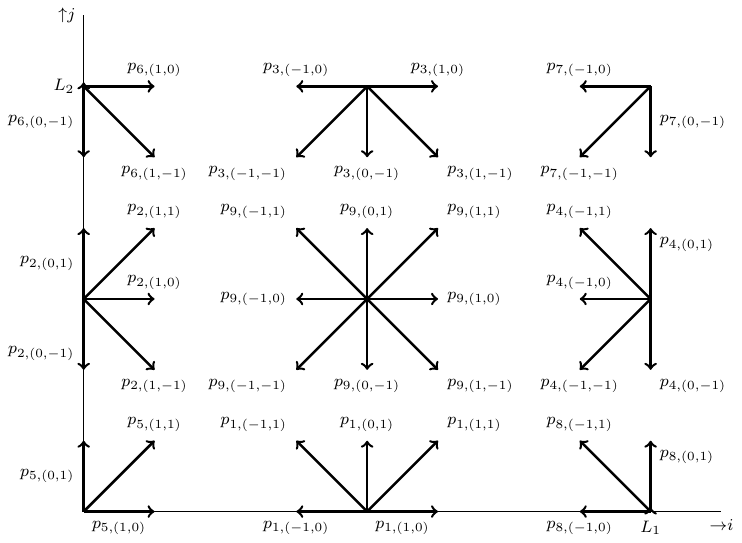}
\hfill{}
\caption{Two-dimensional finite random walk on $S$. The transitions from a state to itself are omitted. \label{fig:rwE}}
\end{figure}

We consider a two-dimensional random walk $R$ on $S$ where 
\begin{equation*}
S = \{0,1,2, \cdots, L_1\} \times\{0,1,2, \cdots, L_2\}.
\end{equation*}
We use a pair of coordinates to represent a state, \ie for $n \in S$, $n = (i,j)$. The state space is naturally partitioned in the following components (see Figure~\ref{fig:partXE}): 
\begin{align*}
C_1 &= \{1,2,3, \cdots, L_1-1\} \times \{0\}, \quad C_2 = \{0\} \times \{1,2,3, \cdots, L_2-1\}, \\
C_3 &= \{1,2,3, \cdots, L_1-1\} \times \{L_2\}, \quad  C_4 = \{L_1\} \times \{1,2,3,\cdots,L_2-1\}, \\ 
C_5 &= \{(0,0)\}, \quad C_6 = \{(0, L_2)\},  \quad C_7 = \{(L_1, L_2)\}, \quad C_8 = \{(L_1, 0)\}, \\
C_9 &= \{1,2,3, \cdots, L_1-1\} \times \{1,2,3, \cdots, L_2-1\}.
\end{align*} 
We refer to this partition as the $C$-partition. The index of the component of state $n \in S$ is denoted by $k(n)$, \ie $n \in C_{k(n)}$. Take for instance, $C_5 = (0,0)$. Then the index of $(0,0)$ is $5$, hence, $k((0,0)) = 5$, \ie $(0,0) \in C_5$. 

Transitions are restricted to the neighboring points (horizontally, vertically and diagonally). For $k = 1,2, \cdots, 9$, we denote by $N_k$ the neighbors of a state in $C_k$. More precisely, $N_1 = \{-1,0,1\} \times \{0,1\}$, $N_2 = \{0,1\} \times \{-1,0,1\}$, $N_3 = \{-1,0,1\} \times \{-1,0\}$, $N_4 = \{-1,0\} \times \{-1,0,1\}$, $N_5 = \{0,1\} \times \{0,1\}$, $N_6 = \{0,1\} \times \{-1,0\}$, $N_7 = \{-1,0\} \times \{-1,0\}$, $N_8 = \{-1,0\} \times \{1,0\}$ and $N_9 = \{-1,0,1\} \times \{-1,0,1\}$. Also, let $N = N_9$.  

Again, let us consider $C_5$. The neighbors, $N_5$, is the product set $\{0,1\} \times \{0,1\}$, which denotes the coordinates of the transitions, either horizontally or vertically.

Let $p_{k,u}$ denote the transition probability from state $n$ in component $k$ to $n + u$, where $u \in N_k$. For $C_5$, we now have $p_{k,u}$ from state $n = (0,0)$ in component $k =5$ to $(0,0) + u$, where $u \in N_5$. This means $u$ could be $(0,0), (0,1), (1,0)$ and $(1,1)$. For instance, $p_{5,(1,0)}$ is the transition probability from state $(0,0)$ in component $5$ to $(0,0) + (1,0)$, \ie (1,0), transition to the right. The transition diagram of a two-dimensional finite random walk can be found in Figure~\ref{fig:rwE}. The transitions from a state to itself are omitted. The system is homogeneous in the sense that the transition probabilities (incoming and outgoing) are translation invariant in each of the components, \ie
\begin{equation}{\label{eq:homogeneousE}}
p_{k(n-u), u} = p_{k(n),u}, \quad \text{for $n - u \in S$ and $u \in N_{k(n)}$}.
\end{equation}
Equation~\eqref{eq:homogeneousE} not only implies that the transition probabilities for each part of the state space are translation invariant but also ensures that also the transition probabilities entering the same component of the state space are translation invariant. 

We assume that the random walk $R$ that we consider is aperiodic, irreducible, positive recurrent, and has invariant probability measure $m(n)$, where $m(n)$ satisfies for all $n \in S$,
\begin{equation*}
m(n) = \sum_{u \in N_{k(n)}} p_{k(n+u), -u} m(n+u).
\end{equation*}

\subsection{Problem formulation}
Our goal is to approximate the steady-state performance of the random walk $R$. The performance measure of interest is
\begin{equation*}
 \mathcal{F} = \sum_{n\in\S} m(n)F(n),
\end{equation*}
where $F(n):\S\to[0,\infty)$  is linear in each of the components from $C$-partition, \ie
\begin{equation}{\label{eq:performancemeasureE}}
F(n) = f_{k(n),0} + f_{k(n),1} i + f_{k(n), 2} j, \quad \text{for} \quad n = (i,j) \in S.
\end{equation}
The constants $f_{k(n),0}$, $f_{k(n),1}$ and $f_{k(n), 2}$ are allowed to be different for different components from the $C$-partition of $S$. 

In general, it is not possible to obtain the probability measure $m(n)$ in a closed-form. Therefore, we will use a perturbed random walk of which the invariant measure has a closed-form expression to approximate the performance measure $\mathcal{F}$.

\begin{figure}
\hfill
\includegraphics{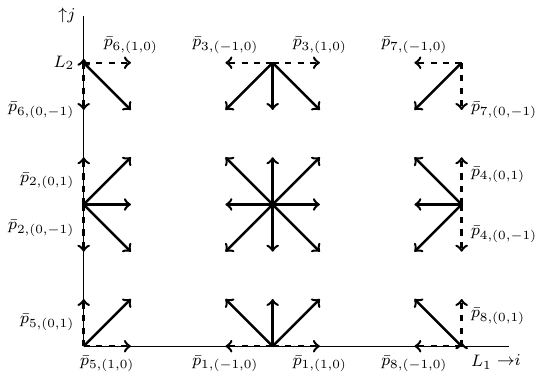}
\hfill{}
\caption{Perturbed random walk $\bar{R}$. \label{fig:rwXE}}
\end{figure}

We approximate the performance measure $\mathcal{F}$ in terms of the perturbed random walk $\bar{R}$. We consider the perturbed random walk $\bar{R}$ in which only the transition probabilities along the boundaries $(C_1, \cdots, C_8)$ are allowed to be different, \ie for instance, $p_{1,(-1,0)}$, $p_{1,(1,0)}$, $p_{1,(0,0)}$ for the state from $C_1$ are allowed to be different in $\bar{R}$, $p_{2,(0,1)}$, $p_{2,(0,-1)}$, $p_{2,(0,0)}$ for the state from $C_2$ are allowed to be different in $\bar{R}$, etc. An example of a perturbed random walk $\bar{R}$ can be found in Figure~\ref{fig:rwXE}. 

We use $\bar{p}_{k,u}$ to denote the probability of $\bar{R}$ jumping from any state $n$ in component $C_k$ to $n + u$, where $u \in N_k$. Moreover, let $q_{k,u} = \bar{p}_{k,u} - p_{k,u}$. The probability measure $\bar{m}$ of $\bar{R}$ is assumed to be of product-form,\ie 
\begin{equation*}
\bar{m}(n) = \alpha \rho^i \sigma^j, 
\end{equation*}
where $n = (i,j)$ for some $(\rho, \sigma) \in (0,1)^2$ and $\alpha \neq 0$. The measure $\bar{m}$ is the invariant measure of $\bar{R}$, \ie it satisfies 
\begin{equation}{\label{eq:balancePE}}
\bar{m}(n) = \sum_{u \in N_{k(n)}} \bar{p}_{k(n+u), -u} \bar{m}(n+u),
\end{equation}
for all $n \in S$.

In the following sections, we are going to find upper and lower bounds of $\mathcal{F}$ in terms of the perturbed random walk $\bar{R}$ defined above.

\section{Proposed approximation scheme}{\label{sec:as4E}}
In this section, we establish an approximation scheme to find upper and lower bounds for performance measures of a two-dimensional finite random walk.

In~\cite{goseling2014linear}, an approximation scheme based on a linear program is developed for a random walk in the quarter-plane. This approximation scheme has also been used in~\cite{chen2015invariant}. We will show in this paper that the technique can be extended to cover our model, \ie a two-dimensional finite random walk. We will explain how this is achieved in the following sections.


\subsection{Markov reward approach to error bounds}

The fact that $R$ and $\bar{R}$ differ only along the boundaries of $S$ makes it possible to obtain the error bounds for the performance measures via the Markov reward approach. An introduction to this technique is provided in~\cite{vandijk11inbook}. We interpret $F$ as a reward function, where $F(n)$ is the one step reward if the random walk is in state $n$. We denote by $F^t(n)$ the expected cumulative reward at time $t$ if the random walk starts from state $n$ at time $0$, \ie

\begin{equation*}
F^t(n) =
\begin{cases}
0, \quad &\text{if } t = 0,\\
F(n) + \sum_{u \in N_{k(n)}} p_{k(n), u} F^{t - 1} (n + u), \quad &\text{if } t > 0,
\end{cases}
\end{equation*}
For convenience, let $F^t(n + u) = 0$ where $u \in \{(s,t)| s,t\in \{-1,0,1\}\}$ if $n + u \notin S$. Terms of the form $F^t(n + u) - F^t(n)$ play a crucial role in the Markov reward approach and are denoted as \emph{bias terms}. Let $D^t_u = F^t(n + u) - F^t(n)$. For the unit vectors $e_1 = (1,0)$, $e_2 = (0,1)$, let $D_1^t(n) = D_{e_1}^t(n)$ and $D_2^t(n) = D_{e_2}^t(n)$.

The next result in~\cite{vandijk11inbook} provides bounds for the approximation error for $\mathcal{F}$. We will use two non-negative functions $\bar{F}$ and $G$ to bound the performance measure $\mathcal{F}$.
\begin{theorem}[~\cite{vandijk11inbook}]{\label{thm:vandijkE}}
Let $\bar{F}$: $S \rightarrow [0, \infty)$ and $G$: $S \rightarrow [0, \infty)$ satisfy
\begin{equation}{\label{eq:requirementE}}
\left|\bar{F}(n) - F(n) + \sum_{u \in N_{k(n)}} q_{k(n), u} D^t_u(n)\right| \leq G(n),
\end{equation}
for all $n \in S$ and $t \geq 0$. Then
\begin{equation}{\label{eq:resultE}}
\sum_{n \in S}[\bar{F}(n) - G(n)] \bar{m}(n) \leq \mathcal{F} \leq \sum_{n \in S}[\bar{F}(n) + G(n)] \bar{m}(n).
\end{equation}
\end{theorem}

\subsection{A linear program approach}
In this section we present a linear program approach to bound the errors. Due to our construction of $\bar{R}$, the random walks $R$ and $\bar{R}$ differ only in the transitions that are along the unit directions, \ie
\begin{equation}{\label{eq:boundaryE}}
q_{k,u} = \bar{p}_{k,u} - p_{k,u} = 0 \quad \text{for} \quad u \neq \{e_1, e_2, -e_1, -e_2, (0,0)\}.
\end{equation}
This restriction will significantly simplify the presentation of the result. 

To start, consider the following optimization problem. We only consider how to obtain the upper bound for $\mathcal{F}$ here because the lower bound for $\mathcal{F}$ can be found similarly.

\noindent \textbf{Problem 1}
\begin{equation}{\label{eq:LPOE}}
\textit{minimize} \quad \sum_{n \in S} [\bar{F}(n) + G(n)] \bar{m}(n),
\end{equation}
\begin{align}
\textit{subject to} \quad &\left|\bar{F}(n) - F(n) + \sum_{s = 1,2} \left(q_{k(n), e_s} D_s^t(n) - q_{k(n), -e_s} D_s^t(n - e_s)\right)\right| \notag \\
                    &\leq G(n), \quad \text{for} \quad n \in S, t \geq 0, \label{eq:LPS1E}\\
                    &\bar{F}(n) \geq 0, G(n) \geq 0, \quad \text{for} \quad n \in S. \label{eq:LPS2E}
\end{align}
The variables in Problem $1$ are the functions $\bar{F}(n)$, $G(n)$ and the parameters are $F(n), \bar{m}(n), q_{k(n), e_s}$ and $D_s^t(n)$ for $n \in S$, $s = 1,2$. Hence, Problem $1$ is a linear programming problem over two non-negative variables $\bar{F}(n)$ and $G(n)$ for every $n \in S$. 

This linear program has infinitely many constraints because we have unbounded time horizon. We will first bound the bias term $D^t_s(n)$ uniformly over $t$. Then we have a linear program with a finite number of variables and constraints. However, further reduction is still needed because the number of variables and constraints will increase rapidly if $L_1$ and $L_2$, which define the size of the state space, increase. Our contribution is to reduce Problem $1$ to a linear programming problem where the number of variables and constraints does not depend on the size of the finite state space. By doing so, we will achieve a constant complexity in the parameters $L_1$ and $L_2$, as opposed to, for instance, the matrix geometric method which has cubic complexity.

We now verify that the objective in Problem $1$ is indeed an upper bound on the performance measure $\mathcal{F}$. Consider $D_{(0,0)}^t(n) = 0$, $D_{-e_s}^t(n) = -D^t_{e_s}(n - e_s)$ for $s = 1,2$ and \eqref{eq:boundaryE}, it follows directly that~\eqref{eq:LPS1E} is equivalent to~\eqref{eq:requirementE}. Therefore, it follows from Theorem~\ref{thm:vandijkE} that the objective of Problem $1$ provides an upper bound on $\mathcal{F}$.

%

\subsection{Bounding the bias terms}
The main difficulty in solving Problem $1$ is the unknown bias terms $D_s^t(n)$. It is in general not possible to find closed-form expressions for the bias terms. Therefore, we introduce two functions $A_s$: $S \rightarrow [0,\infty)$ and $B_s: S \rightarrow [0, \infty)$, $s = 1,2$. We will formulate a finite number of constraints on functions $A_s$ and $B_s$ where $s = 1,2$ such that for any $t$ and $s = 1,2$ we have
\begin{equation}{\label{eq:LPconstrainE}}
-A_s(n) \leq D_s^t(n) \leq B_s(n),
\end{equation}
\ie, the functions $A_s$ and $B_s$ provide bounds on the bias terms uniformly over all $t \geq 0$. In the next section, we will find a finite number of constraints that imply~\eqref{eq:LPconstrainE}. Our method is based on the method that was developed in~\cite{goseling2014linear} for the case of an unbounded state space.

\begin{figure}
\begin{center}
{
\includegraphics{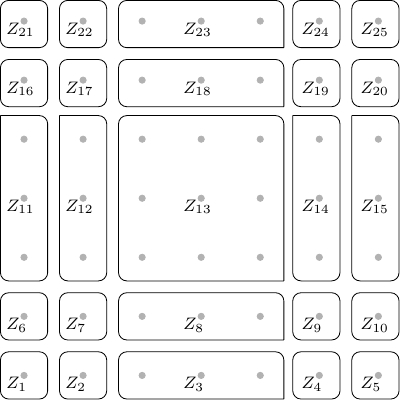}
}
\caption{$Z$-partition of $\S$ with components $Z_1, Z_2, \cdots, Z_{25}$.\label{fig:partitionBE}}
{}
\end{center}
\end{figure}

For notational convenience, as will become clear below, we define a finer partition of $S$, the $Z$-partition. This partition is depicted in Figure~\ref{fig:partitionBE}. For example, we have $Z_1=\{(0,0)\}$, $Z_2=\{(1,0)\}$, $Z_3=\{2, \dots, L_1-2\}\!\!\times\!\!\{0\}$, $Z_4=\{(L_1-1,0)\}$ and $Z_5=\{(L_1,0)\}$, the rest of the elements in the partition are determined similarly. Let $k^z(n)$ denote the label of component from $Z$-partition of state $n \in S$, \ie $n \in Z_{k^z(n)}$. Similar to the definition of $N_k$, let $N^z_k$ denote the neighbors of a state $n$ in $Z_k$ from the $Z$-partition of $S$. 

The constraints which ensure~\eqref{eq:LPconstrainE} are obtained based on an induction in $t$. More precisely, we express $D_s^{t+1}$ as a linear combination of $D_1^t$ and $D_2^t$ as
\begin{equation}{\label{eq:inductionDE}}
D_s^{t+1}(n) = F(n + e_s) - F(n) + \sum_{v = 1,2} \sum_{u \in N^z_{k(n)}} c_{s,k^z(n),v,u} D_v^t(n + u),
\end{equation}
where the $c_{s,k,v,u}$, $s \in \{1,2\}, k \in \{1,2, \cdots, 25\}, v \in \{1,2\}, u \in N^z_k$ are constants. An important property of the $Z$-partition is that starting from any state $n$ in component $k^z$ of the Z-partition the component $k(n+u)$  in the $C$-partition is well defined for all $u\in N^z_k$ and depends only on $k^z$ and $u$. In~\cite{goseling2014linear} it was shown, using this property, that constants $c_{s,k,v,u}$ that ensure~\eqref{eq:inductionDE} always exist and that they can be expressed as simple functions of the transition probabilities of the random walk. The results in~\cite{goseling2014linear} are derived for the random walk on the whole quarter-plane. However, a careful inspection of the results in~\cite{goseling2014linear} reveals that they hold also for our model of a random walk on a bounded state space. Therefore, we refer the reader to~\cite{goseling2014linear} and omit further details here. 

%

We are now ready to bound the bias terms based on~\eqref{eq:inductionDE}. The result, which is easy to verify, states that if $A_s$: $S \rightarrow [0, \infty)$ and $B_s$: $S \rightarrow [0, \infty)$ where $s = 1,2$ satisfy
\begin{multline*}
F(n + e_s) - F(n) \\
+ \sum_{v = 1,2} \sum_{u \in N^z_{k(n)}} \max \{-c_{s,k^z(n),v,u} A_s(n + u), c_{s,k^z(n),v,u} B_s(n + u)\} \leq B_s(n),
\end{multline*}
\begin{multline*}
F(n) - F(n + e_s) \\
+ \sum_{v = 1,2} \sum_{u \in N^z_{k(n)}} \max \{-c_{s,k^z(n),v,u} B_s(n + u), c_{s,k^z(n),v,u} A_s(n + u)\} \leq A_s(n),
\end{multline*}
for all $n \in S$, then
\begin{equation*}{\label{eq:equationCE}}
-A_s (n) \leq D_s^t(n) \leq B_s(n),
\end{equation*}
for $s = 1,2$, $n \in S$ and $t \geq 0$.

%

After bounding the bias terms, we are able to rewrite the linear program Problem $1$ into Problem $2$ with plugging in the upper and lower bounds for $D_s^t(n)$.

\noindent \textbf{Problem 2}
\begin{equation*}{\label{eq:XLPOE}}
\textit{minimize} \quad \sum_{n \in S} [\bar{F}(n) + G(n)] \bar{m}(n),
\end{equation*}
\begin{align*}
\textit{subject to} \quad & \bar{F}(n) - F(n) + \sum_{s = 1,2} \max \{ q_{k(n), e_s} B_s(n) + q_{k(n), -e_s} A_s(n - e_s),  \\
&-q_{k(n), e_s} A_s(n) - q_{k(n), -e_s} B_s(n - e_s) \} \leq G(n), \label{eq:XLPS1E}\\
&F(n) - \bar{F}(n) +  \sum_{s = 1,2} \max \{ q_{k(n), e_s} A_s(n) + q_{k(n), -e_s} B_s(n - e_s), \\
&-q_{k(n), e_s} B_s(n) - q_{k(n), -e_s} A_s(n - e_s) \} \leq G(n) \\
&F(n + e_s) - F(n) + \sum_{v = 1,2} \sum_{u \in N^z_{k(n)}} \max \{-c_{s,k^z(n),v,u} A_s(n + u), \\
&c_{s,k^z(n),v,u} B_s(n + u)\} \leq B_s(n), \\
\displaybreak[4]
&F(n) - F(n + e_s) + \sum_{v = 1,2} \sum_{u \in N^z_{k(n)}} \max \{-c_{s,k^z(n),v,u} B_s(n + u), \\
&c_{s,k^z(n),v,u} A_s(n + u)\} \leq A_s(n), \\
&\bar{F}(n) \geq 0, G(n) \geq 0, A_s(n) \geq 0, B_s(n) \geq 0,\\
&\text{for} \quad n \in S, s \in \{1,2\}.
\end{align*}

\subsection{Fixed number of variables and constraints}

The final step is to reduce Problem $2$ to a linear program with fixed number of variables and constraints regardless of the size of the state space.

We first introduce the notion of a piecewise-linear function on the $Z$-partition. A function $F: S \rightarrow[0, \infty)$ is called $Z$-linear if the function is linear in each of the components from $Z$-partition, \ie
\begin{equation*}
F(n) = f_{k^z(n),0} + f_{k^z(n),1} i + f_{k^z(n), 2} j, \quad \text{for} \quad n = (i,j) \in S.
\end{equation*}
where $f_{k^z(n),0}$, $f_{k^z(n),1}$ and $f_{k^z(n), 2}$ are the constants that define the function. In similar fashion we define $C$-linear functions on the $C$-partition of $S$.

Now, in Problem~$2$ we put the additional constraint that the variables $\bar{F}$, $G$, $A_s$, $B_s$ and $E_s$ are $C$-linear functions. Hence, these functions are defined in terms of variables, the number of which is independent on $L_1$ and $L_2$. Hence, the number of variables in the resulting linear program is independent of $L_1$ and $L_2$.

It remains to show that the number of constraints is independent of $L_1$ and $L_2$. Following the reasoning on the properties of $Z$-partition below~\eqref{eq:inductionDE} it is easy to see that all constraints in Problem~$2$ can be formulated as a non-negativity constraint on a $Z$-linear function. Such a constraint on a $Z$-linear function induces at most $4$ constraints per component in the $Z$-partition, one constraint for each corner of the component. This indicates that the number of constraints does not depend on the size of the state space, since the number of constraints are fixed as well.

\subsection{The optimal solutions}
We are now able to find the upper and lower bounds of $\mathcal{F}$ based on the linear program here.

Let $\mathcal{P}$ denote the set of $(\bar{F}, G)$ for which we are able to find functions $A_s$, $B_s$ and $E_s$ where $s = 1,2$ such that all constraints in Problem $2$ are satisfied. Then, we find the upper and lower bounds for $\mathcal{F}$ as follow.
\begin{equation*}
\mathcal{F}_{up} = \min \left\{\sum_{n \in S} [\bar{F}(n) + G(n)] \bar{m}(n) | (\bar{F}, G) \in \mathcal{P}\right\},
\end{equation*}
and
\begin{equation*}
\mathcal{F}_{low} = \max \left\{\sum_{n \in S} [\bar{F}(n) - G(n)] \bar{m}(n) | (\bar{F}, G) \in \mathcal{P}\right\}.
\end{equation*}
We have now presented the complete approximation scheme to obtain the upper and lower bounds for $\mathcal{F}$ using the perturbed random walk $\bar{R}$ of which the probability measure is of product-form. 

In the following section, we will consider some examples: a tandem queue with finite buffers and some variants of this model.

\section{Application to the Tandem queue with finite buffers}{\label{sec:errorboundsTandemE}}

In this section, we investigate the applications of the approximation scheme proposed in Section~\ref{sec:as4E}.

\subsection{Model description}

Consider a two-node tandem queue with Poisson arrivals at rate $\lambda$. Both nodes have a single server. At most a finite number of jobs, say $L_1$ and $L_2$ jobs, can be present at nodes $1$ and $2$. This includes the jobs in service. An arriving job is rejected if node $1$ is saturated, \ie there are $L_1$ jobs at node $1$. The service time for the jobs at both nodes is exponentially distributed with parameters $\mu_1$ and $\mu_2$, respectively.

\begin{figure}
\begin{center}
\includegraphics{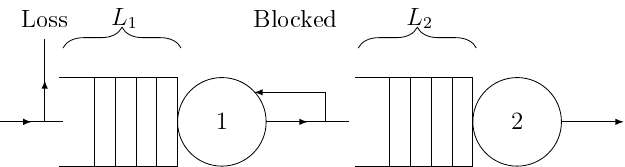}
%
%
%
%
%
%
%
%
\caption{Tandem queue with finite buffers. \label{fig:TandemFBE}}
\end{center}
\end{figure}

When node $2$ is saturated, \ie there are $L_2$ jobs at node $2$, node $1$ stops serving. When it is not blocked, it instantly routes to node $2$. All service times are independent. We also assume that the service discipline is first-in first-out.  

\begin{figure}
\hfill
\includegraphics{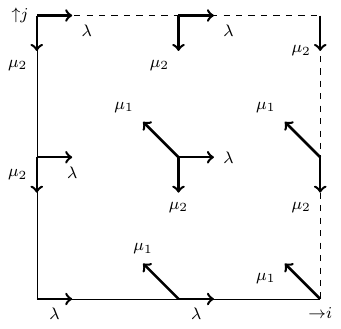}
\hfill{}
\caption{Transition diagram of $R_T$. \label{fig:rwTE}}
\end{figure}

The tandem queue with finite buffers can be represented by a continuous-time Markov process whose state space consists of the pairs $(i,j)$ where $i$ and $j$ are the number of jobs at node $1$ and node $2$, respectively. We now uniformize this continuous-time Markov process to obtain a discrete-time random walk. We assume without loss of generality that $\lambda + \mu_1 + \mu_2 \leq 1$ and uniformize the continuous-time Markov process with uniformization parameter $1$. We denote this random walk by $R_T$. All transition probabilities of $R_T$, except those for the transitions from a state to itself, are illustrated in Figure~\ref{fig:rwTE}.

\subsection{Perturbed random walk of $R_T$}{\label{sec:perturbRWE}}
We now present a perturbed random walk $\bar{R}_T$. The invariant measure of the perturbed random walk $\bar{R}_T$ is of product-form and only the transitions along the boundaries in $\bar{R}_T$ are different from those in $R_T$.

\begin{figure}
\begin{center}
\includegraphics{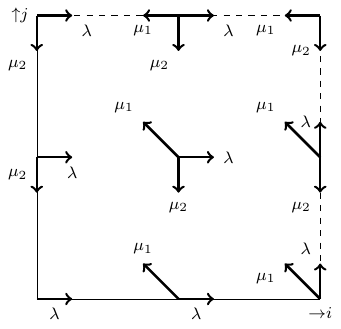}
\end{center}
\caption{Transition diagram of $\bar{R}_T$.\label{fig:PRWE}}
\end{figure}

In the perturbed random walk $\bar{R}_T$, the transition probabilities in the components $C_3, C_4, C_6, C_7, C_8$ are different from those in $R_T$. More precisely, we have $\bar{p}_{3,(1,0)} = \lambda$, $\bar{p}_{3,(-1,0)} = \mu_1$, $\bar{p}_{4,(0,1)} = \lambda$, $\bar{p}_{4,(0,-1)} = \mu_2$, see Figure~\ref{fig:PRWE}. It can be readily verified that the measure, which is of product-form, with $\alpha$, which depends on $L_1$ and $L_2$ as the normalizing constant
\begin{equation*}
\bar{m}(i,j) = \alpha \left(\frac{\lambda}{\mu_1}\right)^i \left(\frac{\lambda}{\mu_2}\right)^j
\end{equation*}
is the probability measure of the perturbed random walk by substitution into the global balance equations~\eqref{eq:balancePE} together with the normalization requirement.

\subsection{Bounding the blocking probability}{\label{sec:examplesE}}
In this section, we provide error bounds for the blocking probability for the tandem queue with finite buffers using our approximation scheme provided in Section~\ref{sec:as4E}. Moreover, we show that our results are better than those obtain by van Dijk et al. in~\cite{vandijk88tandem}. 

For a given performance measure $\mathcal{F}$, we use $\mathcal{F}^{up}$, $\mathcal{F}^{low}$ to denote the upper and lower bounds for $\mathcal{F}$ obtained based on our approximation scheme and $\tilde{\mathcal{F}}^{up}$, $\tilde{\mathcal{F}}^{low}$ to denote the upper and lower bounds based on the method suggested by van Dijk et al.~\cite{vandijk88tandem}.

We use $\mathcal{F}_0$ to denote the blocking probability, \ie the probability that an arriving job is rejected. 
We now consider an example that has also been considered in~\cite{vandijk88tandem}.

\begin{example}{\label{ex:oneE}}
Consider a tandem queue with finite buffers, we have $\lambda = 0.1$, $\mu_1 = 0.2$, $\mu_2 = 0.2$.  
\end{example}

\begin{figure}
  \begin{minipage}{0.25\textwidth}
\begin{tabular}{ |c| }
\hline
  Example~\ref{ex:oneE}  \\
  \hline
  $\lambda = 0.1$  \\
  $\mu_1 = 0.2$  \\
  $\mu_2 = 0.2$  \\
  $L_1 = L_2$ \\
\hline
\end{tabular}    
  \end{minipage}
  \begin{minipage}{0.7\textwidth}
  \includegraphics{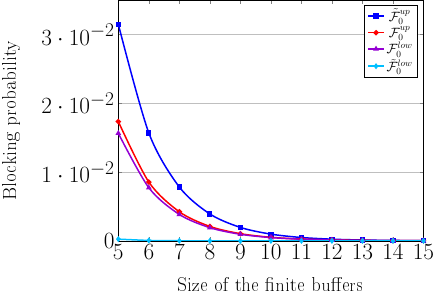}
  \end{minipage}
 \caption{The blocking probability $\mathcal{F}_0$.\label{fig:example1f0E}}  
\end{figure}

\begin{figure}
  \begin{minipage}{0.25\textwidth}
\begin{tabular}{ |c| }
\hline
  Example~\ref{ex:oneE}  \\
  \hline
  $\lambda = 0.1$  \\
  $\mu_1 = 0.2$  \\
  $\mu_2 = 0.2$  \\
  $L_1 = L_2$ \\
\hline
\end{tabular}    
  \end{minipage}
  \begin{minipage}{0.7\textwidth}
  \includegraphics{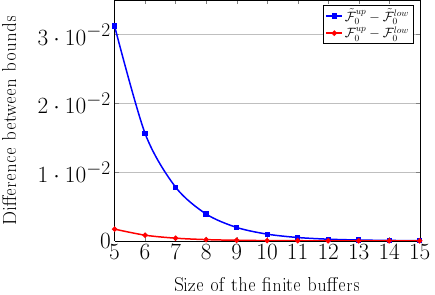}
  \end{minipage}
 \caption{The difference between bounds of $\mathcal{F}_0$.\label{fig:example1f0dE}}  
\end{figure}

We would like to compute the blocking probability of the queueing system. Hence, for the performance measure function $F(n)$, defined in~\eqref{eq:performancemeasureE}, we set the coefficients $f_{k,d}$ where with $k = 1,2,\cdots,9$, $d =0,1,2$ to be $f_{8,0} = 1$, $f_{4,0} = 1$, $f_{7,0} = 1$ and others $0$. The error bounds can be found in Figure~\ref{fig:example1f0E}. Clearly, our results outperform the error bounds obtained in~\cite{vandijk88tandem}. Moreover, the difference between the upper and lower bounds of $\mathcal{F}_0$ are captured in Figure~\ref{fig:example1f0dE}. This indicates that our error bounds are tighter than those in~\cite{vandijk88tandem}.

In addition to the improved bounds, there is another advantage to our method. There is a limitation to the model modification approach that is used in~\cite{vandijk88tandem}. This method requires a different model modification for each specific performance measure. For instance, the specific model modifications which are used to find error bounds for the blocking probability of a tandem queue with finite buffers in~\cite{vandijk88tandem} cannot be used to obtain error bounds for the average number of jobs in the first node. In addition, extra effort is needed to verify that the model modifications are indeed valid for a specific performance measure. In the next section, we will show that our method can easily provide error bounds for other performance measures without extra effort.

\subsection{Bounds for other performance measures}

In this section, we will demonstrate the error bounds for other performance measures for Example~\ref{ex:oneE}, \ie a tandem queue with finite buffers.

Let $\mathcal{F}_1$ be the average number of jobs at node $1$ and $\mathcal{F}_2$ which is the average number of jobs at node $2$.

In general, the models, (\ie the perturbed systems), used to bound the blocking probability in~\cite{vandijk88tandem} cannot be used to bound $\mathcal{F}_1$ and $\mathcal{F}_2$. The method in~\cite{vandijk88tandem} requires different upper and lower bound models for different performance measures. Moreover, this method also requires effort to verify that they are indeed the upper and lower bound models for this specific performance measure. Our approximation scheme does not have this disadvantage. For different performance measure, we only need to change the coefficients $f_{k,d}$ where $k = 1,2,\cdots,9$ and $d =0, 1,2$ in $F(n)$, which is defined in~\eqref{eq:performancemeasureE}.

\begin{figure}
  \begin{minipage}{0.3\textwidth}
\begin{tabular}{ |c| }
\hline
  Example~\ref{ex:oneE}  \\
  \hline
  $\lambda = 0.1$  \\
  $\mu_1 = 0.2$  \\
  $\mu_2 = 0.2$  \\
  $L_1 = L_2$ \\
\hline
\end{tabular}    
  \end{minipage}
  \begin{minipage}{0.65\textwidth}
  \includegraphics{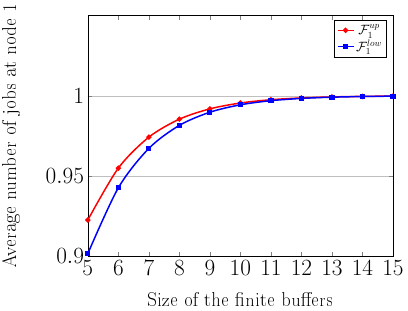}
  \end{minipage}
 \caption{Average number of jobs at node $1$, $\mathcal{F}_1$. \label{fig:example1f1E}}  
\end{figure}

\begin{figure}
  \begin{minipage}{0.3\textwidth}
\begin{tabular}{ |c| }
\hline
  Example~\ref{ex:oneE}  \\
  \hline
  $\lambda = 0.1$  \\
  $\mu_1 = 0.2$  \\
  $\mu_2 = 0.2$  \\
  $L_1 = L_2$ \\
\hline
\end{tabular}    
  \end{minipage}
  \begin{minipage}{0.65\textwidth}
  \includegraphics{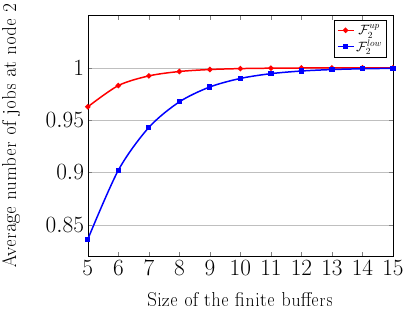}
  \end{minipage}
 \caption{Average number of jobs at node $2$, $\mathcal{F}_2$.{\label{fig:example1f2E}}}  
\end{figure}

It can be readily verified that the performance measure $\mathcal{F}$ is $\mathcal{F}_1$ if and only if we assign following values to the coefficients: $f_{1,1} = 1, f_{8,1} = 1, f_{9,1} = 1, f_{4,1} = 1, f_{3,1} = 1, f_{7,1} = 1$ and others $0$. Figure~\ref{fig:example1f1E} presents the error bounds of $\mathcal{F}_1$.
Similarly, the performance measure $\mathcal{F}$ is $\mathcal{F}_2$ if and only if we assign following values to the coefficients: $f_{2,2} = 1, f_{9,2} = 1, f_{4,2} = 1, f_{6,2} = 1, f_{3,2} = 1, f_{7,2} = 1$ and others $0$. Figure~\ref{fig:example1f2E} presents the error bounds of $\mathcal{F}_2$.

The results show that tight bounds have been achieved with our approximation scheme. Moreover, the only thing we need to change for different performance measures is the input function, which does not require further model modifications. In the next section, we will show that our approximation scheme could also give error bounds for the performance measures of the tandem queue with finite buffers which has a slower or faster server when another node is idle or saturated, respectively, without model modifications as well.

\subsection{Tandem queue with finite buffers and server slow-down/speed-up}

In this section, we consider two variants of the tandem queue with finite buffers. More specifically, we provide error bounds for the blocking probabilities when one server in the tandem with finite buffers is slower or faster if another node is idle or saturated, respectively. 

\subsubsection{Tandem queue with finite buffers and server slow-down}

Tandem queue with server slow-down has been previously studied in, for instance,~\cite{miretskiy2011state,van2005tandem}. A specific type of tandem queue with finite buffers and server slow-down has been considered in~\cite{miretskiy2011state,van2005tandem}. More precisely, the service speed of node $1$ is reduced as soon as the number of jobs in node $2$ reaches some pre-specified threshold because of some sort of protection against frequent overflows.

We consider a different scenario with server slow-down. In our case, the service rate at node $2$ reduces when node $1$ is idle. This comes from a practical situation that when node $1$ is idle, the working pressure for node $2$ decreases and can shift some working capacity to other tasks. Therefore, we consider a two-node tandem queue with Poisson arrivals at rate $\lambda$. Both nodes have a single server. At most a finite number of jobs, say $L_1$ and $L_2$ jobs, can be present at nodes $1$ and $2$, respectively. An arriving job is rejected if node $1$ is saturated. The service time for the jobs at both nodes are exponentially distributed with parameters $\mu_1$ and $\mu_2$, respectively. While node $2$ is saturated, node $1$ stops serving. When it is not blocked, it instantly routes to node $2$. While node $1$ is idle, the service rate of node $2$ becomes $\tilde{\mu}_2$ where $\tilde{\mu}_2 < \mu_2$. All service times are independent. We also assume that the service discipline is first-in first-out.  

\begin{figure}
\hfill
\includegraphics{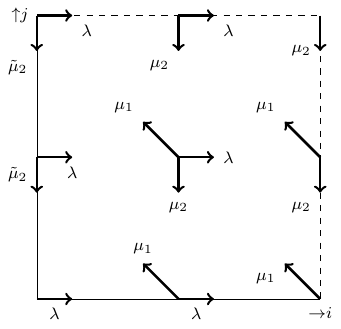}
\hfill{}
\caption{Tandem queue with server slow-down and blocking. \label{fig:rw2E}}
\end{figure}

The tandem queue with finite buffers and server slow-down can be represented by a continuous-time Markov process whose state space consists of the pairs $(i,j)$ where $i$ and $j$ are the number of jobs at node $1$ and node $2$, respectively. We assume without loss of generality that $\lambda + \mu_1 + \mu_2 \leq 1$ and uniformize this continuous-time Markov process with uniformization parameter $1$. Then we obtain a discrete-time random walk. We denote this random walk by $R^{sd}_T$, all transition probabilities of $R^{sd}_T$, except those for the transitions from a state to itself, are illustrated in Figure~\ref{fig:rw2E}. 

It can be readily verified that the random walk $\bar{R}_T$ as defined in Section~\ref{sec:perturbRWE} is a perturbed random walk of $R^{sd}_T$ as well, \ie the transition probabilities in $\bar{R}_T$ only differ from those in $R^{sd}_T$ along the boundaries. We next consider a numerical example.

\begin{example}[slow-down]{\label{ex:twoE}}
Consider a tandem queue with finite buffers and server slow-down, we have $\lambda = 0.1$, $\mu_1 = 0.2$, $\mu_2 = 0.2$ and $\tilde{\mu}_2 = 0.5\mu_2$. 
\end{example}

\begin{figure}
  \begin{minipage}{0.3\textwidth}
\begin{tabular}{ |c| }
\hline
  Example~\ref{ex:twoE}  \\
  \hline
  $\lambda = 0.1$  \\
  $\mu_1 = 0.2$  \\
  $\mu_2 = 0.2$  \\
  $L_1 = L_2$ \\ \hline 
  $\tilde{\mu}_2 = 0.5 \mu_2$ \\
\hline
\end{tabular}    
  \end{minipage}
  \begin{minipage}{0.65\textwidth}
  \includegraphics{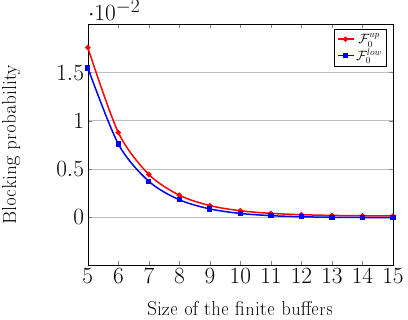}
  \end{minipage}
 \caption{Blocking probability with server slow down.{\label{fig:example2f0E}}}  
\end{figure}

The error bounds for the blocking probability of Example~\ref{ex:twoE} are illustrated in Figure~\ref{fig:example2f0E}.

Notice that our approximation scheme is sufficiently general in the sense that the error bounds for the performance measures of all tandem queue with server slow-down and blocking mentioned in the previous paragraphs can be obtained with our approximation scheme. There are no restrictions on the input random walk.

\subsubsection{Tandem queue with finite buffers and server speed-up}

It is also of great interest to consider a tandem queue with finite buffers and server speed-up.

We consider the following scenario with server speed-up: The service rate at node $2$ increases when node $1$ is saturated. This comes from a practical situation, for instance, when node $1$ is saturated, the working pressure for node $2$ increases to eliminate the jobs in the queueing system. Therefore, we consider a two-node tandem queue with Poisson arrivals at rate $\lambda$. Both nodes have a single server. At most a finite number of jobs, say $L_1$ and $L_2$ jobs, can be present at nodes $1$ and $2$, respectively. An arriving job is rejected if node $1$ is saturated. The service time for the jobs at both nodes are exponential distributed with parameters $\mu_1$ and $\mu_2$, respectively. When node $2$ is saturated, node $1$ stops serving. When it is not blocked, it instantly routes to node $2$. When node $1$ is saturated, the service rate of node $2$ becomes $\bar{\mu}_2$ where $\bar{\mu}_2 > \mu_2$. All service times are independent. We also assume that the service discipline is first-in first-out.  

\begin{figure}
\hfill
\includegraphics{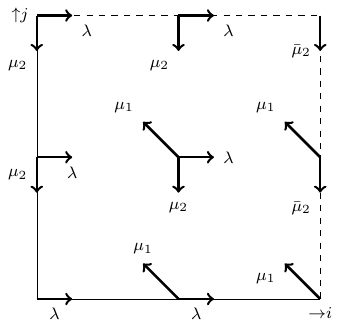}
\hfill{}
\caption{Tandem queue with finite buffers and server speed-up. \label{fig:rw3E}}
\end{figure}

Tandem queue with finite buffers and server speed-up can be represented by a continuous-time Markov process whose state space consists of the pairs $(i,j)$ where $i$ and $j$ are the number of jobs at node $1$ and node $2$, respectively. We assume without loss of generality that $\lambda + \mu_1 + \bar{\mu}_2 \leq 1$ and uniformize this continuous-time Markov process with uniformization parameter $1$. Then we obtain a discrete-time random walk. We denote this random walk by $R^{su}_T$, all transition probabilities of $R^{su}_T$, except those for the transitions from a state to itself, are illustrated in Figure~\ref{fig:rw3E}.

Again, it can be readily verified that the random walk $\bar{R}_T$ as defined in Section~\ref{sec:perturbRWE} is a perturbed random walk of $R^{su}_T$ because only the transitions along the boundaries in $\bar{R}_T$ are different from those in $R^{su}_T$. We next consider the following numerical example.

\begin{example}[speed-up]{\label{ex:threeE}}
Consider a tandem queue with finite buffers and server speed-up, we have $\lambda = 0.1$, $\mu_1 = 0.2$, $\mu_2 = 0.2$ and $\bar{\mu}_2 = 1.2\mu_2$. 
\end{example}

\begin{figure}
  \begin{minipage}{0.3\textwidth}
\begin{tabular}{ |c| }
\hline
  Example~\ref{ex:threeE}  \\
  \hline
  $\lambda = 0.1$  \\
  $\mu_1 = 0.2$  \\
  $\mu_2 = 0.2$  \\
  $L_1 = L_2$ \\ \hline 
  $\bar{\mu}_2 = 1.2 \mu_2$ \\
\hline
\end{tabular}    
  \end{minipage}
  \begin{minipage}{0.65\textwidth}
  \includegraphics{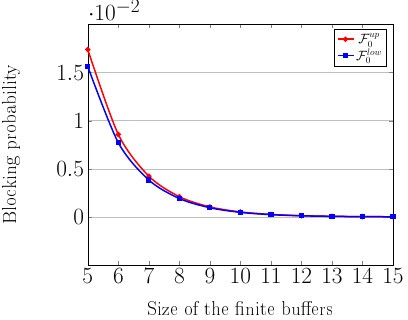}
  \end{minipage}
 \caption{Blocking probability with server speed-up.{\label{fig:example3f0E}}}  
\end{figure}

The error bounds for the blocking probability of Example~\ref{ex:threeE} can be found in Figure~\ref{fig:example3f0E}.

In the next section, we will extend our approximation scheme to the two-dimensional random walk in which one dimension is finite and another dimension is infinite.

\section{Two-node queue with finite buffers at one queue}{\label{sec:onequeue}}
The two-node queue with finite buffers at one queue is a queueing system with two servers, one of them having finite storage capacity. Without loss of generality, we assume node $1$ has finite capacity. If a job arrives at node $1$ when it does not have any more storage capacity, then the job is lost. There is no restriction to the capacity of node $2$.  In general, the two queues influence each other. In particular, the service rate at node $2$ depends on the number of jobs at node $1$. Again we model this queueing system as a two-dimensional random walk for which the state space is finite in one dimension.

We consider a two-dimensional random walk $\tilde{R}$ on $\tilde{S}$ where
\begin{equation*}
\tilde{S} = \{0,1,2, \cdots, L_1\} \times \{0,1,2,3, \cdots\}.
\end{equation*}

Next, we introduce the modified approximation scheme which will be used to find the upper and lower bounds. Similar to the development of the approximation scheme for the two-dimensional finite random walk at both axis, we are able to partition the state space and construct the approximation scheme for the random walk $\tilde{R}$ on state space $\tilde{S}$ based on Markov reward approach. The procedure is different only in the aspect the definition of the components $C_1, C_2, \dots,$ changes. Therefore, we omitted the details and present only the numerical results that have been obtained based on this model.

\section{Application to the coupled-queue with processor sharing and finite buffers at one queue}{\label{sec:errorboundsSharing4E}}

In this section, we apply the approximation scheme to a coupled-queue with processor sharing and finite buffers at one queue. Two coupled processors problem has been extensively studied so far. In particular, Fayolle et al. reduce the problem of finding the generating function of the invariant measure to a Riemann-Hilbert problem in~\cite{fayolle1979two}. However, when we have finite buffers, in general, the methods developed for a coupled-queue with infinite buffers are no longer valid.

\subsection{Model description}

Consider a two-node queue with Poisson arrivals at rate $\lambda_1$ for node $1$ and $\lambda_2$ for node $2$. Both nodes have a single server and at most $L_1$ jobs can be present at nodes $1$ and there is no restriction for the capacity of node $2$. When neither of the nodes is empty they evolve independently, but when one of the queues becomes empty the service rate at another queue changes. An arriving job for node $1$ is rejected when node $1$ is saturated. The service time at both nodes is exponentially distributed with parameters $\mu_1$ and $\mu_2$, respectively, when neither of the queue is empty. When node $1$ is empty, the service rate at node $2$ becomes $\tilde{\mu}_2$ where $\tilde{\mu}_2 > \mu_2$. When node $2$ is empty, the service rate at node $1$ becomes $\tilde{\mu}_1$ where $\tilde{\mu}_1 > \mu_1$. All service requirements are independent. We also assume that the service discipline is first-in first-out.

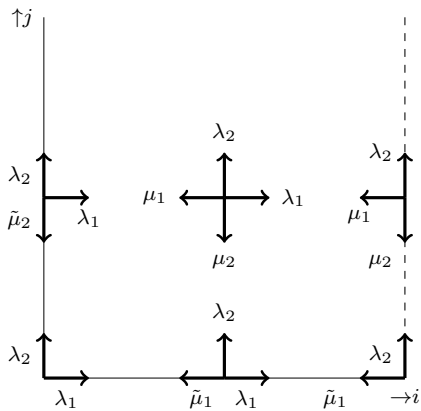
\begin{figure}
\begin{center}
\begin{tikzpicture}[scale=0.6]
\tikzstyle{axes}=[very thin] \tikzstyle{trans}=[very thick,->]
   \draw[axes] (0,0)  -- node[at end, below] {$\scriptstyle \rightarrow i$} (8,0); 
   \draw[axes] (0,0) -- node[at end, left] {$\scriptstyle {\uparrow} {j}$} (0,8);
   \draw[trans] (0,0) to node[below] {$\scriptstyle \lambda_1$} (1,0);
   \draw[trans] (0,0) to node[left] {$\scriptstyle \lambda_2$} (0,1);
   \draw[trans] (4,0) to node[below] {$\scriptstyle \tilde{\mu}_1$} (3,0);
   \draw[trans] (4,0) to node[below] {$\scriptstyle \lambda_1$} (5,0);
   \draw[trans] (4,0) to node[at end, anchor = south]  {$\scriptstyle \lambda_2$} (4,1);
   \draw[trans] (0,4) to node[left] {$\scriptstyle \tilde{\mu}_2$} (0,3);
   \draw[trans] (0,4) to node[at end, anchor = north] {$\scriptstyle \lambda_1$} (1,4);
   \draw[trans] (0,4) to node[left] {$\scriptstyle \lambda_2$} (0,5);
   \draw[trans] (4,4) to node[at end,anchor = west] {$\scriptstyle \lambda_1$} (5,4);
   \draw[trans] (4,4) to node[at end, anchor = south] {$\scriptstyle \lambda_2$} (4,5);
   \draw[trans] (4,4) to node[at end, anchor = east] {$\scriptstyle \mu_1$} (3,4);
   \draw[trans] (4,4) to node[at end, anchor = north] {$\scriptstyle \mu_2$} (4,3);
   \draw[dashed] (8,8) -- (8,0);
   \draw[trans] (8,4) to node[at end, anchor = north east] {$\scriptstyle \mu_2$}(8,3);
   \draw[trans] (8,4) to node[at end, anchor = north] {$\scriptstyle \mu_1$}(7,4);
   \draw[trans] (8,4) to node[at end, anchor = east] {$\scriptstyle \lambda_2$}(8,5);
   \draw[trans] (8,0) to node[at end, anchor = north east] {$\scriptstyle \tilde{\mu}_1$}(7,0);
   \draw[trans] (8,0) to node[at end, anchor = north east] {$\scriptstyle \lambda_2$}(8,1);
   
\end{tikzpicture} 
\end{center}
\caption{Coupled-queue with processor sharing and finite buffers at one queue.\label{fig:rw4E}}
\end{figure}

This coupled-queue with processor sharing and finite buffers at one queue can be represented by a continuous-time Markov process whose state space consists of the pairs $(i,j)$ where $i$ and $j$ are the number of jobs at node $1$ and node $2$, respectively. We assume without loss of generality that $\lambda_1 + \lambda_2 + \tilde{\mu}_1 + \tilde{\mu}_2 \leq 1$ and uniformize this continuous-time Markov process with uniformization parameter $1$. Then we obtain a discrete-time random walk. We denote this random walk by $R_C$. All transition probabilities of $R_C$, except those for the transitions from a state to itself, are illustrated in Figure~\ref{fig:rw4E}.

\subsection{Perturbed random walk $\bar{R}_C$}
We now display a perturbed random walk $\bar{R}_C$ of $R_C$ such that the probability measure of $\bar{R}_C$ is of product-form and only the transitions along the boundaries in $\bar{R}_C$ are different from those in $R_C$.

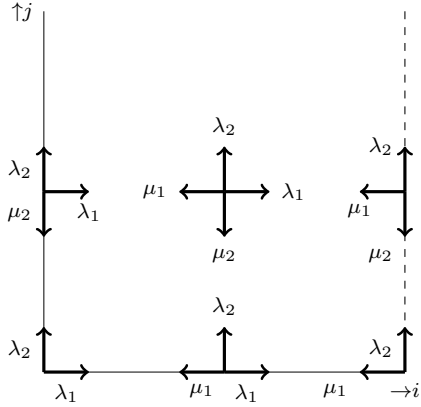
\begin{figure}
\begin{center}
\begin{tikzpicture}[scale=0.6]
\tikzstyle{axes}=[very thin] \tikzstyle{trans}=[very thick,->]
   \draw[axes] (0,0)  -- node[at end, below] {$\scriptstyle \rightarrow i$} (8,0); 
   \draw[axes] (0,0) -- node[at end, left] {$\scriptstyle {\uparrow} {j}$} (0,8);
   \draw[trans] (0,0) to node[below] {$\scriptstyle \lambda_1$} (1,0);
   \draw[trans] (0,0) to node[left] {$\scriptstyle \lambda_2$} (0,1);
   \draw[trans] (4,0) to node[below] {$\scriptstyle \mu_1$} (3,0);
   \draw[trans] (4,0) to node[below] {$\scriptstyle \lambda_1$} (5,0);
   \draw[trans] (4,0) to node[at end, anchor = south]  {$\scriptstyle \lambda_2$} (4,1);
   \draw[trans] (0,4) to node[left] {$\scriptstyle \mu_2$} (0,3);
   \draw[trans] (0,4) to node[at end, anchor = north] {$\scriptstyle \lambda_1$} (1,4);
   \draw[trans] (0,4) to node[left] {$\scriptstyle \lambda_2$} (0,5);
   \draw[trans] (4,4) to node[at end,anchor = west] {$\scriptstyle \lambda_1$} (5,4);
   \draw[trans] (4,4) to node[at end, anchor = south] {$\scriptstyle \lambda_2$} (4,5);
   \draw[trans] (4,4) to node[at end, anchor = east] {$\scriptstyle \mu_1$} (3,4);
   \draw[trans] (4,4) to node[at end, anchor = north] {$\scriptstyle \mu_2$} (4,3);
   \draw[dashed] (8,8) -- (8,0);
   \draw[trans] (8,4) to node[at end, anchor = north east] {$\scriptstyle \mu_2$}(8,3);
   \draw[trans] (8,4) to node[at end, anchor = north] {$\scriptstyle \mu_1$}(7,4);
   \draw[trans] (8,4) to node[at end, anchor = east] {$\scriptstyle \lambda_2$}(8,5);
   \draw[trans] (8,0) to node[at end, anchor = north east] {$\scriptstyle \mu_1$}(7,0);
   \draw[trans] (8,0) to node[at end, anchor = north east] {$\scriptstyle \lambda_2$}(8,1);
   
\end{tikzpicture} 
\end{center}
\caption{Transition diagram of the perturbed random walk $\bar{R}_C$.\label{fig:rw4PE}}
\end{figure}

It can be readily verified that the invariant measure of the perturbed random walk $\bar{R}_C$ in Figure~\ref{fig:rw4PE},
which is of product-form, with $\alpha$, which depends on $L_1$ as the normalizing constant
\begin{equation*}
\bar{m}(n) = \alpha \left(\frac{\lambda_1}{\mu_1}\right)^i \left(\frac{\lambda_2}{\mu_2}\right)^j \quad \text{where} \quad n = (i,j),
\end{equation*}
is the probability measure of the perturbed random walk by substitution into the global balance equations~\eqref{eq:balancePE} together with the normalization requirement.

We next illustrate a numerical example of a coupled-queue with processor sharing and finite buffers at one queue.

\subsection{Numerical results}

\begin{example}{\label{ex:fourE}}
Consider a coupled-queue with finite buffers at one queue, we have $\lambda_1 = \lambda_2 = 0.15, \mu_1 = \mu_2 = 0.2$, $\tilde{\mu}_1 = \tilde{\mu}_2 = 0.25$.
\end{example}

%

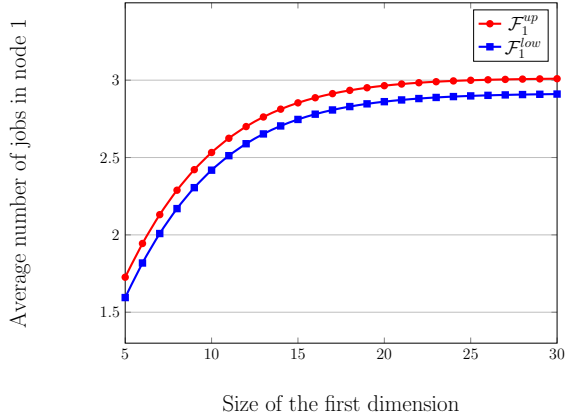
\begin{figure}
  \begin{minipage}{0.3\textwidth}
\begin{tabular}{ |c| }
\hline
  Example~\ref{ex:fourE}  \\
  \hline
  $\lambda_1 = 0.15$  \\
  $\lambda_2 = 0.15$ \\
  $\mu_1 = 0.2$  \\
  $\mu_2 = 0.2$  \\ \hline 
  $\tilde{\mu}_1 = 1.25 \mu_1$ \\
  $\tilde{\mu}_2 = 1.25 \mu_2$ \\
\hline
\end{tabular}    
  \end{minipage}
  \begin{minipage}{0.65\textwidth}
%
%
%
%
\begin{tikzpicture}[scale = 0.5]

\definecolor{mycolor1}{rgb}{0.5804,0,0.8275}
\definecolor{mycolor2}{rgb}{0,0.7490,1.0000}

\begin{axis}[%
view={0}{90},
width=4.52083333333333in,
height=3.565625in,
scale only axis,
xmin=5, xmax=30,
ymin=1.3, ymax=3.5,
x label style={at={(axis description cs:0.5,-0.08)},anchor=north},
y label style={at={(axis description cs:-0.1,.5)},anchor=south},
xtick = {5,10,15,20,25,30},
xlabel={$\scriptstyle \text{\Large Size of the first dimension}$},
ytick = {1.5,2,2.5,3},
ylabel={$\scriptstyle \text{\Large Average number of jobs in node $1$}$},
ymajorgrids]

\addplot [
color=red,
ultra thick,
solid,
mark=*,
mark options={solid},
smooth
]
coordinates{
 (5,1.72632961308)(6,1.94430220402)(7,2.13070015681)(8,2.28869839921)(9,2.42154812905)(10,2.53242488151)(11,2.62432917932)(12,2.70002488015)(13,2.76200550614)(14,2.81248156448)(15,2.85338326119)(16,2.88637395361)(17,2.91287051023)(18,2.93406752757)(19,2.9509630815)(20,2.96438433738)(21,2.97501188788)(22,2.98340212446)(23,2.99000728063)(24,2.99519302654)(25,2.99925365852)(26,3.00242503197)(27,3.00489544474)(28,3.00681470514)(29,3.00830162255)(30,3.00945014926)
};
\addlegendentry{\large $\mathcal{F}_1^{up}$}

\addplot [
color=blue,
ultra thick,
solid,
mark=square*,
mark options={solid},
smooth
]
coordinates{


 (5,1.5947278828)(6,1.81814990663)(7,2.00861039289)(8,2.16977861734)(9,2.30519202832)(10,2.41819930497)(11,2.51191096249)(12,2.58916407618)(13,2.65250303932)(14,2.70417509378)(15,2.74613789371) (16,2.78007592363)(17,2.80742273044)(18,2.82938635851)(19,2.84697591719)(20,2.86102775339)(21,2.87223018817)(22,2.88114617845)(23,2.88823357718)(24,2.89386289373)(25,2.89833261301)(26,2.90188223152)(27,2.90470322419)(28,2.90694818067)(29,2.90873835264)(30,2.91016984296)
};
\addlegendentry{\large $\mathcal{F}_1^{low}$}

\end{axis}
\end{tikzpicture}
  \end{minipage}
 \caption{Average number of jobs in node $1$.{\label{fig:example_new_2}}}  
\end{figure}

We approximate the average number of jobs in node $1$. We use $F_1$ to denote the average number of jobs in node $1$. The upper and lower bounds of $F_1$, which are denoted by $F_{1}^{up}$ and $F_{1} ^{low}$, can be found in Figure~\ref{fig:example_new_2}.

We see from the results in Figure~\ref{fig:example_new_2} that our approximation scheme can also be extended to finite random walks at one axis. Moreover, note that when $L_1$, \ie the size of the first dimension, is increasing, the values of the upper and lower bounds reach a limit.

In the next numerical example, we will fix the service rate. We present the error bounds for the corresponding performance measure when the occupation rate, \ie $\rho = \frac{\lambda}{\mu}$ increases, even close to $1$.

\begin{example}{\label{ex:fiveE}}
Consider a coupled-queue with finite buffers at one queue, we have $\mu_1 = \mu_2 = 0.2$, $\tilde{\mu}_1 = \tilde{\mu}_2 = 0.25$, $L_1=20$. Let $\rho$ changes from $0.5$ to $0.95$.
\end{example}

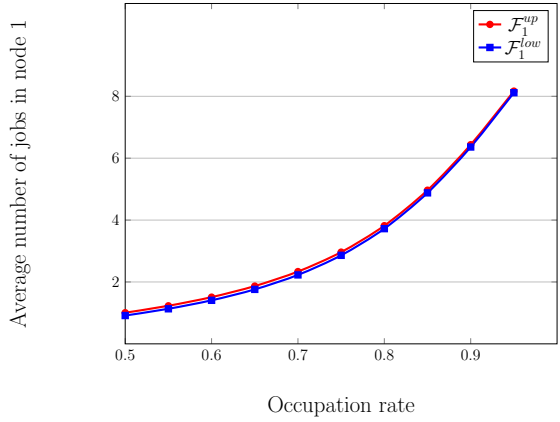
\begin{figure}
  \begin{minipage}{0.3\textwidth}
\begin{tabular}{ |c| }
\hline
  Example~\ref{ex:fiveE}  \\
  \hline
  $\rho = 0.5 \dots 0.95$ \\
  $\mu_1 = 0.2$  \\
  $\mu_2 = 0.2$  \\ \hline 
  $\tilde{\mu}_1 = 1.25 \mu_1$ \\
  $\tilde{\mu}_2 = 1.25 \mu_2$ \\ \hline
  $L_1 = 20$ \\
\hline
\end{tabular}    
  \end{minipage}
  \begin{minipage}{0.65\textwidth}
%
%
%
%
\begin{tikzpicture}[scale = 0.5]

\definecolor{mycolor1}{rgb}{0.5804,0,0.8275}
\definecolor{mycolor2}{rgb}{0,0.7490,1.0000}

\begin{axis}[%
view={0}{90},
width=4.52083333333333in,
height=3.565625in,
scale only axis,
xmin=0.5, xmax=1,
ymin=0, ymax=11,
x label style={at={(axis description cs:0.5,-0.08)},anchor=north},
y label style={at={(axis description cs:-0.1,.5)},anchor=south},
xtick = {0.5,0.6,0.7,0.8,0.9},
xlabel={$\scriptstyle \text{\Large Occupation rate}$},
ytick = {2,4,6,8},
ylabel={$\scriptstyle \text{\Large Average number of jobs in node $1$}$},
ymajorgrids]

\addplot [
color=red,
ultra thick,
solid,
mark=*,
mark options={solid},
smooth
]
coordinates{
 (0.5,1.00628913426)(0.55,1.230138989)(0.6,1.50958975766)(0.65,1.8672738531)(0.7,2.33605594356)(0.75,2.96438433738)(0.8,3.81888002961)(0.85,4.96577624584)(0.9,6.43579818573)(0.95,8.16744298981)};
\addlegendentry{\large $\mathcal{F}_1^{up}$}

\addplot [
color=blue,
ultra thick,
solid,
mark=square*,
mark options={solid},
smooth
]
coordinates{

 (0.5,0.915147672746)(0.55,1.13297110557)(0.6,1.40705101912)(0.65,1.76008894756)(0.7,2.22805569482)(0.75,2.86102775339)(0.8,3.72287063376)(0.85,4.87975938225)(0.9,6.36122412175)(0.95,8.11780325411)};
\addlegendentry{\large $\mathcal{F}_1^{low}$}

\end{axis}
\end{tikzpicture}
  \end{minipage}
 \caption{Average number of jobs in node $1$ when $\rho$ increases.{\label{fig:example_new_3}}}  
\end{figure}

We see from Figure~\ref{fig:example_new_3} that the error bounds are quite tight as well.

Next, we present several examples for blocking probabilities, which is again denoted by $F_0$, based on Example~\ref{ex:fiveE} in which the size of the buffers in the first dimension increases from $20$ to $10000$.

\begin{example}{\label{ex:sixE}}
Consider a coupled-queue with finite buffers at one queue, we have $\mu_1 = \mu_2 = 0.2$, $\tilde{\mu}_1 = \tilde{\mu}_2 = 0.25$, $L_1 = 20$ and the occupation rate increases from $0.5$ to $0.95$.
\end{example}

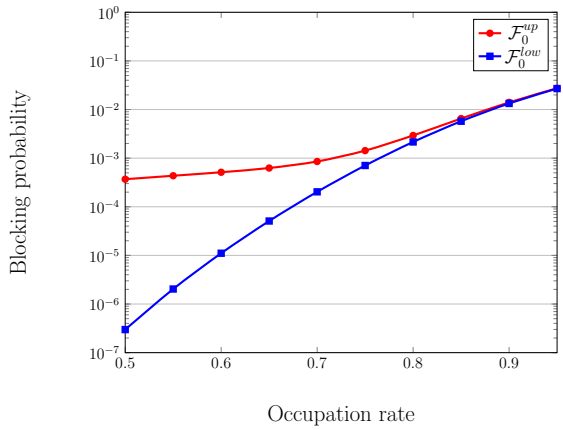
\begin{figure}
  \begin{minipage}{0.3\textwidth}
\begin{tabular}{ |c| }
\hline
  Example~\ref{ex:sixE}  \\
  \hline
  $\rho = 0.5 \dots 0.95$ \\
  $\mu_1 = 0.2$  \\
  $\mu_2 = 0.2$  \\ \hline 
  $\tilde{\mu}_1 = 1.25 \mu_1$ \\
  $\tilde{\mu}_2 = 1.25 \mu_2$ \\ \hline
  $L_1 = 20$ \\
\hline
\end{tabular}    
  \end{minipage}
  \begin{minipage}{0.65\textwidth}
%
%
%
%
\begin{tikzpicture}[scale = 0.5]

\definecolor{mycolor1}{rgb}{0.5804,0,0.8275}
\definecolor{mycolor2}{rgb}{0,0.7490,1.0000}

\begin{semilogyaxis}[%
view={0}{90},
width=4.52083333333333in,
height=3.565625in,
scale only axis,
xmin=0.5, xmax=0.95,
ymin=1e-7, ymax=1,
x label style={at={(axis description cs:0.5,-0.08)},anchor=north},
y label style={at={(axis description cs:-0.1,.5)},anchor=south},
xtick = {0.5,0.6,0.7,0.8,0.9},
xlabel={$\scriptstyle \text{\Large Occupation rate}$},
ylabel={$\scriptstyle \text{\Large Blocking probability}$},
ymajorgrids]

\addplot [
color=red,
ultra thick,
solid,
mark=*,
mark options={solid},
smooth
]
table[
	x index=0, y index=2
	]
	{blocking_pr_rho_L20.csv};

\addlegendentry{\large $\mathcal{F}_0^{up}$}

\addplot [
color=blue,
ultra thick,
solid,
mark=square*,
mark options={solid},
smooth
]
table[
	x index=0, y index=1
	]
	{blocking_pr_rho_L20.csv};

\addlegendentry{\large $\mathcal{F}_0^{low}$}

\end{semilogyaxis}
\end{tikzpicture}
  \end{minipage}
\caption{The converted blocking probability ($y = \log Y$), $L_1 = 20$.{\label{fig:blocking_pr_rho_L20_log}}}  
\end{figure}

The bounds for blocking probabilities are very close in this case, hence, we convert these probabilities by applying logarithm to the $y$ axis in Figure~\ref{fig:blocking_pr_rho_L20_log} and also in following examples.

\begin{example}{\label{ex:sevenE}}
Consider a coupled-queue with finite buffers at one queue, we have $\mu_1 = \mu_2 = 0.2$, $\tilde{\mu}_1 = \tilde{\mu}_2 = 0.25$, $L_1 = 500$ and the occupation rate increases from $0.98$ to $0.99$.
\end{example}



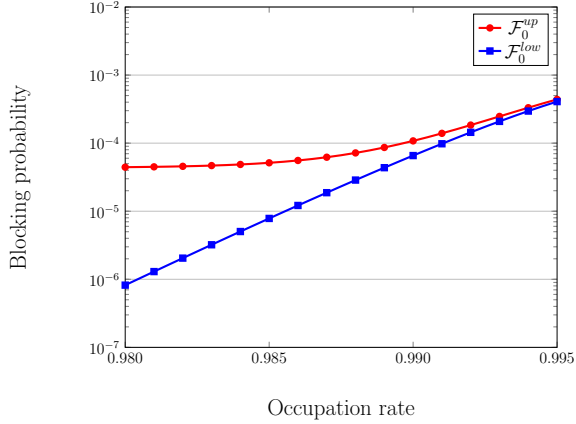
\begin{figure}
  \begin{minipage}{0.3\textwidth}
\begin{tabular}{ |c| }
\hline
  Example~\ref{ex:sevenE}  \\
  \hline
  $\rho = 0.980 \dots 0.995$ \\
  $\mu_1 = 0.2$  \\
  $\mu_2 = 0.2$  \\ \hline 
  $\tilde{\mu}_1 = 1.25 \mu_1$ \\
  $\tilde{\mu}_2 = 1.25 \mu_2$ \\ \hline
  $L_1 = 500$ \\
\hline
\end{tabular}    
  \end{minipage}
  \begin{minipage}{0.65\textwidth}
%
%
%
%
\begin{tikzpicture}[scale = 0.5]

\definecolor{mycolor1}{rgb}{0.5804,0,0.8275}
\definecolor{mycolor2}{rgb}{0,0.7490,1.0000}

\begin{semilogyaxis}[%
x tick label style={
	    /pgf/number format,
	    fixed,
	    fixed zerofill,
	    precision=3
	},
view={0}{90},
width=4.52083333333333in,
height=3.565625in,
scale only axis,
xmin=0.980, xmax=0.995,
ymin=1e-7, ymax=1e-2,
x label style={at={(axis description cs:0.5,-0.08)},anchor=north},
y label style={at={(axis description cs:-0.1,.5)},anchor=south},
xtick = {0.980,0.985,0.990,0.995},
xlabel={$\scriptstyle \text{\Large Occupation rate}$},
ylabel={$\scriptstyle \text{\Large Blocking probability}$},
ymajorgrids]

\addplot [
color=red,
ultra thick,
solid,
mark=*,
mark options={solid},
smooth
]
table[
	x index=0, y index=2
	]
	{blocking_pr_rho_L500.csv};

\addlegendentry{\large $\mathcal{F}_0^{up}$}

\addplot [
color=blue,
ultra thick,
solid,
mark=square*,
mark options={solid},
smooth
]
table[
	x index=0, y index=1
	]
	{blocking_pr_rho_L500.csv};

\addlegendentry{\large $\mathcal{F}_0^{low}$}

\end{semilogyaxis}
\end{tikzpicture}
  \end{minipage}
\caption{The converted blocking probability ($y = \log Y$), $L_1 = 500$.{\label{fig:blocking_pr_rho_L500_loglog}}}  
\end{figure}

Next, we also extend these numerical results to the case when $L_1 = 10000$.

\begin{example}{\label{ex:eightE}}
Consider a coupled-queue with finite buffers at one queue, we have $\mu_1 = \mu_2 = 0.2$, $\tilde{\mu}_1 = \tilde{\mu}_2 = 0.25$, $L_1 = 10000$ and the occupation rate increases from $0.98$ to $0.99$.
\end{example}



\begin{figure}
  \begin{minipage}{0.3\textwidth}
\begin{tabular}{ |c| }
\hline
  Example~\ref{ex:eightE}  \\
  \hline
  $\rho = 0.9990 \dots 0.9998$ \\
  $\mu_1 = 0.2$  \\
  $\mu_2 = 0.2$  \\ \hline 
  $\tilde{\mu}_1 = 1.25 \mu_1$ \\
  $\tilde{\mu}_2 = 1.25 \mu_2$ \\ \hline
  $L_1 = 10000$ \\
\hline
\end{tabular}    
  \end{minipage}
  \begin{minipage}{0.65\textwidth}
%
%
%
%
\begin{tikzpicture}[scale = 0.5]

\definecolor{mycolor1}{rgb}{0.5804,0,0.8275}
\definecolor{mycolor2}{rgb}{0,0.7490,1.0000}

\begin{semilogyaxis}[%
x tick label style={
	    /pgf/number format,
	    fixed,
	    fixed zerofill,
	    precision=4
	},
view={0}{90},
width=4.52083333333333in,
height=3.565625in,
scale only axis,
xmin=0.9990, xmax=0.9998,
ymin=1e-8, ymax=1e-3,
x label style={at={(axis description cs:0.5,-0.08)},anchor=north},
y label style={at={(axis description cs:-0.1,.5)},anchor=south},
xtick = {0.9990,0.9992,0.9994,0.9996,0.9998,1},
xlabel={$\scriptstyle \text{\Large Occupation rate}$},
ylabel={$\scriptstyle \text{\Large Blocking probability}$},
ymajorgrids]

\addplot [
color=red,
ultra thick,
solid,
mark=*,
mark options={solid},
smooth
]
table[
	x index=0, y index=2
	]
	{blocking_pr_rho_L10000.csv};

\addlegendentry{\large $\mathcal{F}_0^{up}$}

\addplot [
color=blue,
ultra thick,
solid,
mark=square*,
mark options={solid},
smooth
]
table[
	x index=0, y index=1
	]
	{blocking_pr_rho_L10000.csv};

\addlegendentry{\large $\mathcal{F}_0^{low}$}

\end{semilogyaxis}
\end{tikzpicture}
  \end{minipage}
\caption{The converted blocking probability ($y = \log Y$), $L_1 = 10000$.{\label{fig:blocking_pr_rho_L10000_loglog}}}  
\end{figure}
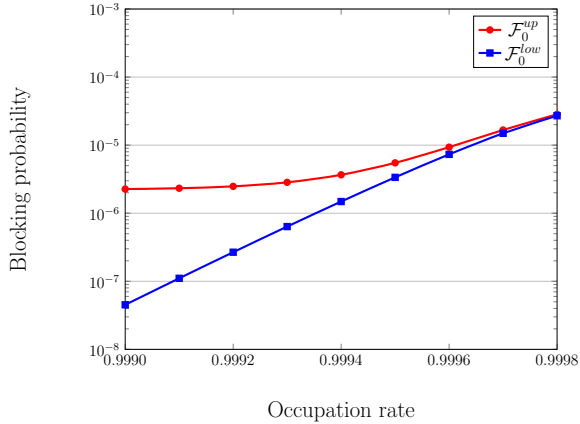

We see from the above examples that relatively tight bounds are obtained efficiently based on our approach. As discussed in the introduction that the matrix geometric method has cubic complexity in $L_1$.

\section{Conclusion}{\label{sec:conclusionE}}

In this paper, we presented a general approximation scheme for a two-node queue with finite buffers at either one or both queues, which establishes error bounds for a large class of performance measures. Our work is an extension of the linear programming approach developed in~\cite{goseling2014linear} to approximate performance measures of random walks in the quarter-plane. 

We first developed an approximation scheme for a two-node queue with finite buffers at both queues. We then applied this approximation scheme to obtain bounds for performance measures of a tandem queue in which both buffers are finite and some variants of this model. We also extended the approximation scheme to deal with a two-node queue with finite buffers at only one queue. We applied our approximation scheme to a coupled-queue with finite buffers at one queue. The approximation scheme gives tight bounds for various performance measures, like the blocking probability and the average number of jobs at node $1$. We also obtain error bounds for the blocking probabilities when the size of the buffers in one dimension is really large.

To summarize, the complexity for solving a system of linear equations is at least $O(L_1^2)$ and the variations of matrix geometric method share a complexity of $O(L_1^3)$. Therefore, when $L_1$ is large, our approach, of which the complexity is a constant in $L_1$, acts as a promising alternative to finding the invariant measures.

\section{Acknowledgment}
Yanting Chen acknowledges support through the NSFC grant 71701066, the Fundamental Research Funds for the Central Universities and a CSC scholarship [No. 2008613008]. Xinwei Bai acknowledges support through a CSC scholarship [No. 201407720012]. This work is partly supported by the Netherlands Organization for Scientific Research (NWO) grant $612.001.107$.





%
%
%

\bibliographystyle{plain}
\bibliography{bibfile} 

\end{document}